\documentclass[12 pt]{amsart}
\usepackage{amsfonts,amsmath,amsthm,amssymb, amscd}
\usepackage[top=2cm, bottom=5cm, left=3cm, right=1cm]{geometry}

\usepackage{tikz-cd}
\usepackage{xcolor}
\usepackage{placeins}
\usepackage{hyperref}
\usepackage[all]{xy}

\newtheorem{theorem}{Theorem}[section]
\theoremstyle{definition}

\newtheorem{corollary}[theorem]{Corollary}
\newtheorem{lemma}[theorem]{Lemma}
\newtheorem{proposition}[theorem]{Proposition}

\newtheorem{definition}[theorem]{Definition}
\newtheorem{example}[theorem]{Example}

\usepackage{lineno}
\usepackage{color}

\title[Orders of Products of Slanted Class Transpositions]
{
Orders of Products of Slanted Class Transpositions}
\author[V.~G.~Bardakov, A.~L.~Iskra]{V.~G.~Bardakov, A.~L.~Iskra}

\address{Sobolev Institute of Mathematics, 4 Acad. Koptyug avenue, 630090, Novosibirsk, Russia.}
\address{Novosibirsk State Agrarian University, Dobrolyubova street, 160, Novosibirsk, 630039, Russia.}
\address{Regional Scientific and Educational Mathematical Center of Tomsk State University, 36 Lenin Ave., Tomsk, Russia.}
\email{bardakov@math.nsc.ru}

\address{Sobolev Institute of Mathematics, 4 Acad. Koptyug avenue, 630090, Novosibirsk, Russia.}
\email{a.iskra@g.nsu.ru}
\date{\today}

\begin{document}
\maketitle
\begin{abstract}
In the present work, we continue the research initiated in the preprint: V. G. Bardakov, A. L. Iskra, Orders of products of horizontal class transpositions, arXiv:2409.13341, and related to S. Kohl’s question on the orders of products of pairs of class transpositions. In the preprint, an answer was given to S. Kohl’s question for horizontal class transpositions. In the present work, the products of pairs of slanted class transpositions are considered, and under certain conditions, their orders are determined, with it being established that the number of different orders is finite.


\end{abstract}
\maketitle
\tableofcontents

\section{Introduction}\label{s1}

S. Kohl \cite{K} introduced the group $CT(\mathbb{Z})$, which is a subgroup of the permutation group $Sym(\mathbb{Z})$ of the set of integers $\mathbb{Z}$. The group $CT(\mathbb{Z})$ is generated by class transpositions, it is a simple group, it contains all symmetric groups
 $S_n$
and hence all finite groups. A class transposition is an involution that interchanges residue classes modulo $m_1$ and $m_2$.  
More precisely, for a natural number $m > 1$, for every $r$ satisfying the inequalities $0 \leq r < m$, the residue class modulo $m$ is defined as  
$$
r(m) = r + m \mathbb{Z} = \{ r + km~|~k \in \mathbb{Z} \}.
$$
If $r_1(m_1) \cap r_2(m_2) = \emptyset$, the class transposition $\tau_{r_1(m_1),r_2(m_2)}$ is defined as the permutation which interchanges $r_1+km_1$ and $r_2+km_2$ for each integer  $k$ and which fixes all other numbers.  
It is easy to verify that the permutation $\tau_{r_1(m_1), r_2(m_2)}$ is a class transposition if and only if the greatest common divisor $\gcd(m_1, m_2)$ does not divide $|r_1 - r_2|$.  

The class transposition $\tau_{r_1(m_1), r_2(m_2)}$ is the product of independent transpositions:  
$$
\tau_{r_1(m_1),r_2(m_2)} = \prod_{k \in \mathbb{Z}} (r_1 +  m_1 k, r_2 +  m_2 k).
$$

In the Kourovka notebook \cite{Kour}, S. Kohl was posed the question (see Question 18.48): is it true that there exists only a finite set of natural numbers that occur as the orders of the product of two class transpositions? Using computer calculations, it was established that the product of a pair of class transpositions can have an order belonging to the following set 
$$
\{ 1, 2, 3, 4, 6, 8, 10, 12, 15, 20, 24, 30, 40, 42, 60, 84, 120, 168, 420 \}.
$$
All these numbers divide 840. S. Kohl then formulated the question \cite{K2}: is this true for an any pair of class transpositions whose product has finite order?

In the work \cite{B}, each class transposition $\tau = \tau_{r_1(m_1),r_2(m_2)}$ was associated with a segment $AB$ lying on the plane and connecting the point $A = (r_1, m_1)$ with the point $B = (r_2, m_2)$. The points $A$ and $B$ are called the vertices of the class transposition $\tau$. The permutation $\tau_{r_1(m_1), r_2(m_2)}$ is called horizontal if $m_1 = m_2$. In this case, the segment corresponding to the class transposition is horizontal (parallel to the $OX$ axis). A class transposition that is not horizontal is called {\it slanted}.  
Horizontal class transpositions generate the subgroup $CT_{\infty}$, introduced in \cite{K}. The group of generalized class transpositions, which contains the group $CT(\mathbb{Z})$, was introduced in \cite{NN}.

In the work \cite{B}, it was proved that the order of the product of a pair of horizontal class transpositions belongs to the set $\{1,2,3,4,6,12\}$, and for every number in this set there exists a pair of horizontal class transpositions whose product has that order.

We denote the order of an element $g$ of a group $G$ by $ord(g)$. In the present work, the orders of the products of pairs of slanted class transpositions are studied. In particular, the products of pairs of slanted class transpositions that share a common vertex are examined. For these, it has been proved the following theorem 

\begin{theorem}
Let $\tau_{r(m),r_1(m_1)}$ and $\tau_{r(m),r_2(m_2)}$ be class transpositions with a common vertex. Then $ord(\tau_{r(m),r_1(m_1)}\cdot\tau_{r(m),r_2(m_2)})\in\{1,3,\infty\}$, and moreover $ord(\tau_{r(m),r_1(m_1)}\cdot\tau_{r(m),r_2(m_2)})=\infty$ if and only if
$$
r_1(m_1) \cap r_2(m_2) \neq \emptyset~\mbox{and}~r_1(m_1)\neq r_2(m_2).
$$
\end{theorem}

Also studied are pairs of class transpositions for which the projections of the segments onto the $OX$ or $OY$ axis coincide. In other words, in the first case, pairs of class transpositions $\tau_{l(m_1),r(m_2)}$ and $\tau_{l(m_3),r(m_4)}$ are considered, which we shall call  {\it equal-residue} class transpositions and in the second case, pairs of class transpositions $\tau_{r_1(m),r_2(n)}$ and $\tau_{r_3(m),r_4(n)}$ are considered, which we shall call  {\it equal-modulus} class transpositions. For equal-residue and equal-modulus class transpositions, the following theorem has been proved.

\begin{theorem}
1) Let $\tau_{l(m_1),r(m_2)}$ and $\tau_{l(m_3),r(m_4)}$ be equal-residue class transpositions and assume that $\frac{m_1}{m_2}\neq \frac{m_3}{m_4}$. Then $ord(\tau_{l(m_1),r(m_2)}\cdot \tau_{l(m_3),r(m_4)})=\infty$.

2) Let $\tau_{r_1(m),r_2(n)}$ and $\tau_{r_3(m),r_4(n)}$ be equal-modulus class transpositions and assume that $r_1\leqslant r_4,\,r_3\leqslant r_2$. Then $ord(\tau_{r_1(m),r_2(n)}\cdot\tau_{r_3(m),r_4(n)})\in \{1,2,3,6,\infty\}$.
\end{theorem}

The question of describing the orders of the products of equal-residue and equal-modulus class transpositions that do not satisfy the conditions of Theorem 1.2 remains open.

The structure of the work is as follows. In  \S~\ref{s2}, we recall the construction of the graph determined by a pair of class transpositions that carries information about the order of their product. In \S~\ref{s3}, Theorem 1.1 (see Theorem  \ref{3.1}) is proved, and a criterion is provided for when the product of a pair of class transpositions with a common vertex has infinite order. In \S~\ref{s4}, Theorem 1.2 (see Theorem  \ref{4.3}) is proved.

\bigskip

%%%%%%%%%%%%%%%%%%%%%%%%%%%%%%%%%%

\section{The graph of the product of two class transpositions}\label{s2}

In the work \cite{B} a graph corresponding to the product of a pair of class transpositions was introduced, and a connection was established between the connected components of this graph and the orbits of the product. We recall the necessary definitions. Consider a pair of class transpositions
$$
\tau_1 = \tau_{r_1(m_1),r_2(m_2)} =  \prod\limits_{k \in \mathbb{Z}} (r_1+m_1k, r_2+m_2k),~~~\tau_2 =\tau_{r_3(m_3),r_4(m_4)} = \prod\limits_{l \in \mathbb{Z}} (r_3+m_3l , r_4+m_4 l).
$$
To construct the graph, choose a pair of disjoint sets
$$
V_1 = \{ a_k, b_k ~| ~ k \in \mathbb{Z}\},~~~V_2 = \{ c_l, d_l ~| ~ l \in \mathbb{Z}\}.
$$
That is, the set $V_i$ contains two disjoint subsets, each of which is indexed by the integers. Define the graph $\Gamma(\tau_1,\tau_2)=(V, E)$, where the vertex set $V=V_1\sqcup V_2$ is the disjoint union of $V_1$ and $V_2$, which are associated with $\tau_1$ and $\tau_2$ via the map 
 $$
 \mu \colon V \to Supp(\tau_1) \cup Supp(\tau_2) \subseteq \mathbb{Z},
 $$
defined by the rule
$$
\mu(a_k) = r_1+m_1k,~~\mu(b_k) = r_2+m_2k,~~\mu(c_l) = r_3+m_3 l,~~\mu(d_l) = r_4+m_4 l.
$$
Evidently, the restriction of $\mu$ to each $V_i$, $i = 1, 2$, is a bijection.
The set of undirected edges $E$ consists of pairs $\{v_i, v_j \}\in E$, for which one of the following three cases holds:

1) $ \mu(v_i) =  \mu(v_j)$ and $v_i\in V_1, v_j\in V_2$ or $v_i\in V_2, v_j\in V_1$,

2)$ \tau_1(\mu(v_i))= \mu(v_j)$ and $v_i,v_j\in V_1$,

3) $\tau_2(\mu(v_i))= \mu(v_j)$ and $v_i, v_j\in V_2$.

Edges satisfying condition (1) are called {\it edges of the first type}, while edges satisfying conditions (2) or (3) are called {\it edges of the second type}. 
For brevity, we will write $\tau_i(u) = v$  instead of $\tau_i(\mu(u)) = \mu(v)$, where $u, v \in V$ and $i = 1, 2$.
If $\Delta$ is a subgraph of $\Gamma(\tau_1,\tau_2)$,  we define
$$
\mu(\Delta)=\{\mu(h)\,|\,h\text{ is a vertex of } \Delta \}.
$$

The following lemma is evident.

\begin{lemma}
For any vertex $v$  of the graph $\Gamma(\tau_1,\tau_2)$ one of the following two statements holds: 

1) The vertex $v$  has degree 1 and, in this case, is incident only to an edge of the second type;

2) The vertex $v$ has degree 2 and, in this case, is incident to an edge of the first type and an edge of the second type.
\end{lemma}

Using this lemma, it is easy to classify the connected components of the graph  $\Gamma(\tau_1,\tau_2)$. If a connected component is finite (i.e. contains a finite number of vertices) and all its vertices have degree 2, then choose an arbitrary vertex and denote it by $v_1$. In this case, the connected component is represented by a sequence of vertices
$$
v_1 v_2\dots v_{n-1}v_1,
$$
in which every pair of consecutive vertices is connected by an edge. Thus, if a connected component is given by a sequence of vertices, we assume that consecutive vertices are connected by an edge. Such connected components are called connected components of the {\it first type}. If a finite connected component contains vertices of degree 1, then it is easy to see that there are exactly two such vertices. One is denoted by $v_1$ and the other by $v_n$  so that the inequality $\mu(v_1)<\mu(v_n)$ holds. Such connected components are written in the form 
$$
v_1v_2\dots v_n.
$$
and are called connected components of the {\it second type}. If all vertices of a connected component have degree 2 and there are infinitely many of them, then, after choosing any vertex of this component and denoting it by $v_0$  we write the connected component in the form 
$$
\prod_{i\in\mathbb{Z}}v_i.
$$
Such connected components are called connected components of the {\it third type}. If a connected component contains an infinite number of vertices and there exists a vertex of degree 1, then that vertex is unique. Denote it by $v_1$ and write the connected component in the form
$$
\prod_{i=1}^{\infty}v_i.
$$
Connected components of this kind are called connected components of the {\it fourth type}. The {\it length} of a connected component is defined as the number of edges it contains.

In the work \cite{B} the following was established:

\begin{theorem} \label{t1}
There exists a mapping $\psi$ from the set of connected components $\{C_i\}_{i\in I}$ of the graph $\Gamma(\tau_1,\tau_2)$ into $Sym(\mathbb{Z})$ such that

1)
\begin{equation*}
\psi(C_i) = 
 \begin{cases}
   (u_1\ldots u_{\frac{n}{4}})(v_1\ldots v_{\frac{n}{4}}), &\text{\it{if $C_i$ is of the first type of length n}},\\
   (w_1\ldots w_{\frac{n+3}{2}}), &\text{\it{if $C_i$  is of the second type of length n}},\\
   (\ldots,u_{-1},u_0,u_1,\ldots)(\ldots,v_{-1},v_0,v_1,\ldots), &\text{\it{if $C_i$ is of the third type}},\\
   (\ldots w_{-1},w_0,w_1,\ldots), &\text{\it{if $C_i$ is of fourth type}}.
 \end{cases}
\end{equation*}

2) $Supp(\psi(C_i))=\mu(C_i)$.

3) $\prod_{i\in I}\psi(C_i)=\tau_1\cdot\tau_2.$

\end{theorem}

\bigskip

%%%%%%%%%%%%%%%%%%%%%%%%%%%%%%%%%%

\section{The product of two class transpositions with a common vertex}\label{s3}

We will say that the class transpositions $\tau_{r_1(m_1),r_2(m_2)}$ and $\tau_{r_3(m_3),r_4(m_4)}$ have a {\it common vertex} if one of the following conditions holds:

1) $r_1=r_3,\,m_1=m_3$;

2) $r_1=r_4,\,m_1=m_4$;

3) $r_2=r_3,\,m_2=m_3$;

4) $r_2=r_4,\,m_2=m_4$.

Geometrically, this means that the intersection of the segments corresponding to these class transpositions either consists of a single vertex or the segments coincide.
Taking into account that $\tau_{r_i(m_i),r_j(m_j)} = \tau_{r_j(m_j),r_i(m_i)}$, when considering a pair of transpositions with a common vertex, we will assume that condition 1) holds.

\begin{theorem} \label{3.1}
Let $\tau_{r(m),r_2(m_2)}$ and $\tau_{r(m),r_4(m_4)}$ be  class transpositions with a common vertex. Then $ord(\tau_{r(m),r_2(m_2)}\cdot\tau_{r(m),r_4(m_4)})\in\{1,3,\infty\}$ and $ord(\tau_{r(m),r_2(m_2)}\cdot\tau_{r(m),r_4(m_4)})=\infty$ if and only if 
$$
r_2(m_2) \cap r_4(m_4) \neq \emptyset~\mbox{and}~r_2(m_2)\neq r_4(m_4).
$$
\end{theorem}

\begin{proof}

{\it Case 1:} $r_2(m_2)\cap r_4(m_4)=\emptyset$. Fix $k\in \mathbb{Z}$. By the definition of the graph, the vertices $a_k$ and $b_k$ as well as the vertices $c_k$ and $d_k$ are connected by edges. By the condition $\mu(a_k)=\mu(c_k)$, then the vertices $a_k$ and $c_k$ are connected by an edge. There are no other edges among the vertices $a_k, b_k, c_k$ and $d_k$. It is easy to see that no edges emanate from these vertices to any other vertices of the graph. Consequently, the graph $\Gamma(\tau_{r(m),r_1(m_1)},\tau_{r(m),r_2(m_2)})$ contains the connected component $d_k c_k a_k b_k$. Since  $k$ is an arbitrary integer, we obtain
$$
\Gamma(\tau_{r(m),r_1(m_1)},\tau_{r(m),r_2(m_2)})=\bigsqcup_{k\in\mathbb{Z}} d_k c_k a_k b_k.
$$ 
By Theorem \ref{t1}, we have $ord(\tau_{r(m),r_2(m_2)}\cdot\tau_{r(m),r_4(m_4)})=3$.

{\it Case 2:} $r_2(m_2)=r_4(m_4)$, i.e. the class transpositions coincide and hence 
$$
ord(\tau_{r(m),r_2(m_2)}\cdot\tau_{r(m),r_4(m_4)})=ord(\tau_{r(m),r_2(m_2)}\cdot\tau_{r(m),r_2(m_2)})=1.
$$

{\it Case 3:} $r_2(m_2)\cap r_4(m_4)\neq\emptyset$ and $r_2(m_2)\neq r_4(m_4)$.
As in Case 1, for every integer~$k$ we have $\{a_k,b_k\}\in E$, $\{c_k,d_k\}\in E$, and $\{a_k,c_k\}\in E$; that is, the graph will contain subgraphs of the form $d_k c_k a_k b_k$.  Let us show that  the graph $\Gamma(\tau_{r(m),r_1(m_1)},\tau_{r(m),r_2(m_2)})$ contains a connected component of length $\geq n$. To do it it is enough to show that the following system 
\begin{equation*}
 \begin{cases}
   b_{x_1}=d_{x_2}, 
   \\
   b_{x_2}=d_{x_3},
   \\
   \cdots\cdots\cdots
   \\
   b_{x_{n}}=d_{x_{n+1}},
 \end{cases}
\end{equation*}
which by   the implications
$$
b_{x_i}=d_{x_{i+1}}\Leftrightarrow r_2+m_2x_i=r_4+m_4 x_{i+1}  \Leftrightarrow m_2 x_i - m_4 x_{i+1} = r_4 - r_2,
$$
 is equivalent to the following system 
of Diophantine equations   
\begin{equation*}
 \begin{cases}
   m_2x_1-m_4x_2=r_4-r_2, 
   \\
   m_2x_2-m_4x_3=r_4-r_2,
   \\
   \cdots\cdots\cdots\cdots\cdots\cdots\cdot\cdot
   \\
   m_2x_n-m_4x_{n+1}=r_4-r_2,
 \end{cases}
\end{equation*}
has a solution $ (x_{1}, \ldots, x_{n+1}) = (x_{1,n}(t_n),...,x_{n+1,n}(t_n))$, for which

1) $x_{n+1,n}(t_n)=\tilde{x}_{n+1}+m_2^nt_n$;

2) there is a set $T_n \subseteq \mathbb{Z}$
such that for all $t_n \in \mathbb{Z}\setminus T_n$ the following inequalities hold $x_{i,n}(t_n)\neq x_{j,n}(t_n)$ if $i\neq j$.

Induction by $n$. For $n=1$ we have one equation
$$
 m_2 x_1 - m_4 x_2 = r_4 - r_2.
$$
Since the residue classes $r_2(m_2)$ and $r_4(m_4)$ intersect, this equation is solvable, i.e., $gcd(m_2,m_4)$ divides $|r_4 - r_2|$.
 Without loss of generality, assume that $gcd(m_2,m_4)=1$. 
The pair
$$
x_{1,1}(t_1)=\tilde{x}_1+m_4t_1,\,x_{2,1}(t_1)=\tilde{x}_2+m_2t_1
$$
is a solution of the equation $b_{x_1}=d_{x_2}$ for every $t_1\in\mathbb{Z}$, provided that $(\tilde{x}_1,\tilde{x}_2)$ is a particular solution of the equation. Consider the set
$$
T_1=\{\tilde{t}_1\in\mathbb{Z}\,|\,x_{1,1}(\tilde{t}_1)=x_{2,1}(\tilde{t}_1)\}.
$$
Note that $|T_1|\leq 1$. For any $t_1\in\mathbb{Z}\setminus T_1$ we have $x_{1,1}(t_1)\neq x_{2,1}(t_1)$. Consequently, the graph $\Gamma(\tau_{r(m),r_2(m_2)},\tau_{r(m),r_4(m_4)})$ contains subgraphs of the form
$$
d_{x_{1,1}(t_1)}c_{x_{1,1}(t_1)}a_{x_{1,1}(t_1)}b_{x_{1,1}(t_1)}d_{x_{2,1}(t_1)}c_{x_{2,1}(t_1)}a_{x_{2,1}(t_1)}b_{x_{2,1}(t_1)}.
$$

Suppose further that the following system of $n-1$ equations
\begin{equation*}
 \begin{cases}
   m_2x_1-m_4x_2=r_4-r_2, 
   \\
   m_2x_2-m_4x_3=r_4-r_2,
   \\
   \cdots\cdots\cdots\cdots\cdots\cdots\cdot\cdot
   \\
   m_2x_{n-1}-m_4x_n=r_4-r_2.
 \end{cases}
\end{equation*}
have a solution $(x_{1,n-1}(t_{n-1}),\ldots,x_{n,n-1}(t_{n-1})) $, where $t_{n-1}\in\mathbb{Z}\setminus T_{n-1},\,|T_{n-1}|\leq 1$ and $x_{n,{n-1}}(t_{n-1})=\tilde{x}_n+m_2^{n-1}t_{n-1}$, and such that $x_{i,{n-1}}(t_{n-1})\neq x_{j,{n-1}}(t_{n-1})$ for all $i\neq j$ and all $t_{n-1}\in\mathbb{Z}\setminus T_{n-1}$. 

We will show that the system
\begin{equation*}
 \begin{cases}
   m_2x_1-m_4x_2=r_4-r_2, 
   \\
   m_2x_2-m_4x_3=r_4-r_2,
   \\
   \cdots\cdots\cdots\cdots\cdots\cdots\cdot\cdot
   \\
   m_2x_n-m_4x_{n+1}=r_4-r_2
 \end{cases}
\end{equation*}
has a solution $x_{1,n}(t_n),...,x_{n+1,n}(t_n)$, where $t_n\in\mathbb{Z}\setminus T_n$ and $x_{n+1,n}(t_n)=\tilde{x}_{n+1}+m_2^nt_n$, such that $x_{i,n}(t_n)\neq x_{j,n}(t_n)$ for all $i\neq j$ and all $t_n\in\mathbb{Z}\setminus T_n$.

Substitute in the last equation $x_n=x_{n,n-1}(t_{n-1})=\tilde{x}_n+m_2^{n-1}t_{n-1}$; we obtain 
$$
m_2^nt_{n-1}-m_4x_{n+1}=r_4-r_2-m_2^{n-1}\tilde{x}_n.
$$
Since $gcd(m_2^n,m_4)=1$, a solution exists and is given by
$$
t_{n-1}(t_n)=\tilde{t}_{n-1}+m_4t_n,\,x_{n+1,n}(t_n)=\tilde{x}_{n+1}+m_2^nt_n.
$$
Denote
$$
T_n=\{\tilde{t}_n\in\mathbb{Z}\,|\,t_{n-1}(\tilde{t}_n)\in T_{n-1}.\}
$$
It is easy to see that $|T_n|\leq 1.$ Define 
$$ 
x_{1,n}(t_n)=x_{1,n-1}(t_{n-1}(t_n)),\ldots,\,x_{n,n}(t_n)=x_{n,n-1}(t_{n-1}(t_n)),\,x_{n+1,n}(t_n)=\tilde{x}_{n+1}+m_2^nt_n,
$$
where $t_n\in\mathbb{Z}\setminus T_n$. This is the desired solution. 

We will show that $x_{i,n}(t_n)\neq x_{j,n}(t_n)$ for all $i\neq j$ and all $t_n\in\mathbb{Z}\setminus T_n$. If $i,j<n+1$, then  
$$
x_{i,n}(t_n)=x_{i,n-1}(t_{n-1}(t_n))\neq x_{j,n-1}(t_{n-1}(t_n))=x_{j,n}(t_n).
$$
This inequality follows from the fact that $t_{n-1}(t_n)\in \mathbb{Z}\setminus T_{n-1}$ for every $t_n\in\mathbb{Z}\setminus T_n$.

Now, let $i = n + 1$. Suppose that there exists some $j\geq 1$ such that $x_{n+1,n}(t_n)=x_{j,n}(t_n)$ for some $t_n\in \mathbb{Z}\setminus T_n.$ If $j>1$, then 
$$
r_2+m_2x_{j-1,n}(t_n)=r_4+m_4x_{j,n}(t_n)=r_4+m_4x_{n+1,n}(t_n)=r_2+m_2x_{n,n}(t_n),
$$
but $x_{j-1,n}(t_n)\neq x_{n,n}(t_n)$. This is a contradiction.

Now, let $j=1$. Define  
$$
\tilde{b}_l = r_2+m_2x_{l,n}(t_n),~~\tilde{d}_s = r_4+m_4x_{s,n}(t_n).
$$
Note that $\tilde{b}_l=\tilde{d}_{l+1}$. Moreover, $\tilde{b}_l<\tilde{b}_s$ if and only if $\tilde{d}_l<\tilde{d}_s$. Since $x_{1,n}(t_n)\neq x_{2,n}(t_n)$ it follows that $\tilde{b}_1\neq \tilde{b}_2$. Suppose that $\tilde{b}_1<\tilde{b}_2$. Then $\tilde{d}_2<\tilde{d}_3$, which implies that $\tilde{b}_2<\tilde{b}_3$. Continuing this process, we obtain the inequalities 
$$\tilde{b}_1<\tilde{b}_2<\ldots<\tilde{b}_n,~~\tilde{d}_1<\tilde{d}_2<\ldots<\tilde{d}_{n+1}.
$$ 
Then 
$$
\tilde{b}_1<\tilde{b}_n=\tilde{d}_{n+1}=\tilde{d}_1<\tilde{d}_2=\tilde{b}_1,
$$
a contradiction.

Thus, we have shown that for every natural number $n$ the above system of Diophantine equations has a solution. It follows that for every natural number $n$ the graph $\Gamma(\tau_{r(m),r_2(m_2)},\tau_{r(m),r_4(m_4)})$   contains a subgraph of the form
$$
\prod_{i=1}^nb_{x_{i,n}(t_n)}c_{x_{i,n}(t_n)}a_{x_{i,n}(t_n)}d_{x_{i,n}(t_n)}.
$$
By Theorem \ref{t1}, the order of the product $\tau_{r(m),r_2(m_2)}\cdot\tau_{r(m),r_4(m_4)}$ is infinite.
\end{proof}

From the proven theorem the following corollary follows:

\begin{corollary}
Let $\tau_{r(m),r_2(m_2)}$ and $\tau_{r(m),r_4(m_4)}$ be class transpositions with a common vertex. Then the order of their product is equal to 3 if and only if the intersection $r_2(m_2) \cap r_4(m_4)$ is empty.
\end{corollary}

\bigskip

The theorem proven does not provide information about the cycles contained in the product of two class transpositions with a common vertex. The following proposition gives a partial answer to this question. To denote an infinite periodic sequence, we introduce the following notation:
$$
 {\bf P}_{l=s}^{\infty} \langle f_1(l), f_2(l), \ldots f_n(l) \rangle =  f_1(s), f_2(s), \ldots f_n(s), f_1(s+1), f_2(s+1), \ldots f_n(s+1), \ldots.
$$

\begin{proposition}
Let $k$ be an odd, $m = 2 m_1$ be even, and $m\geq 4$. Then $\tau_{0(2),k(m)}$ is a class transposition, and the product $\tau_{0(2),1(2)}\cdot\tau_{0(2),k(m)}$ acts on non-negative integers by the rules 
$$
 (\ldots, 2,  3, {\bf P}_{t=1}^{\infty} \langle k + (k-1) \frac{m_1^{t-1} - 1}{m_1 - 1}m_1 + 2 m_1^t\rangle )~\mbox{if}~ k \not= 3,
$$
$$
(\ldots, 2, 0, 1, {\bf P}_{t=1}^{\infty} \langle 3 + 2 \frac{m_1^{t-1} - 1}{m_1 - 1}m_1\rangle)~\mbox{if}~ k = 3
$$
and, in particular, it has infinite order.
\end{proposition}

\begin{proof}
Depending on  the value of $k$, we consider two cases.	

{\it Case 1}: $k=1$. First, suppose that $m = 2m_1$  is even. In this case, the product $\tau_{0(2),1(2)}\cdot\tau_{0(2),k(m)}$ contains the following infinite cycle
$$
(\ldots, 2, 3, 1+2m_1, 1+2m_1^2, 1+2m_1^3, \ldots ) =  (\ldots, 2,  3, {\bf P}_{t=1}^{\infty} \langle 1 + 2 m_1^t\rangle ).
$$
Now, suppose that $m = 1 + 2m_1$ is odd.  In this case, since $gcd(2, m) = 1$, the permutation $\tau_{0(2),k(m)}$ is not a class transposition. 

{\it Case 2}: $k>1$. For the same reasons as above, the odd $m$ case is impossible, so we must have $m = 2m_1$. If $k$ were also even, then this case would not occur; therefore, we assume that $k$ is odd. First, suppose that $k=3$. Then the product $\tau_{0(2),1(2)}\cdot\tau_{0(2),k(m)}$ contains the cycle:
$$
(\ldots, 2, 0, 1, 3, 3+ 2 m_1, 3 + 2(1 + m_1)m_1, 3 + 2(1 +  m_1 + m_1^2) m_1, \ldots ) = (2, 0, 1, {\bf P}_{t=1}^{\infty} \langle 3 + 2 \frac{m_1^{t-1} - 1}{m_1 - 1}m_1\rangle).
$$
Now, suppose that $k>3$. Then the product $\tau_{0(2),1(2)}\cdot\tau_{0(2),k(m)}$ contains the cycle:
$$
(\ldots, 2, 3, k+ 2 m_1, k + (k-1 + 2m_1)m_1, k + (k-1 + ((k-1) + 2 m_1)m_1)m_1, \ldots ) = 
$$
$$
= (\ldots, 2,  3, {\bf P}_{t=1}^{\infty} \langle k + (k-1) \frac{m_1^{t-1} - 1}{m_1 - 1}m_1 + 2 m_1^t\rangle ).
$$
\end{proof}

\section{The product of two equal-residue and equal-modulus class transpositions} \label{s4}

We introduce the following

\begin{definition}
1) Class transpositions $\tau_{r_1(m_1),r_2(m_2)}$ and $\tau_{r_3(m_3),r_4(m_4)}$ are called {\it equal-residue} if $\{ r_1, r_2 \} = \{ r_3, r_4\}$.

2) Class transpositions $\tau_{r_1(m_1),r_2(m_2)}$ and $\tau_{r_3(m_3),r_4(m_4)}$ are called {\it equal-modulus} if $\{ m_1, m_2 \} = \{ m_3, m_4\}$. 
\end{definition}

Furthermore, when considering equal-residue (equal-modulus) class transpositions, without loss of generality we assume that $r_1 = r_3$ and $r_2 = r_4$ (respectively, $m_1 = m_3$ and $m_2 = m_4$).
Note that the set of all class transpositions is partitioned into disjoint equivalence classes, where two class transpositions are equivalent if they are equal-residue. Similarly, the set of all class transpositions is partitioned into disjoint equivalence classes, where two class transpositions are equivalent if they are equal-modulus.

Denote by $T$ the set of class transpositions, and on the set of pairs $(\tau_1, \tau_2) \in T \times T$   we introduce the following predicates. We say that

$(\tau_1, \tau_2) \in K_g$, if $\tau_1$ and $\tau_2$ are horizontal class transpositions;

$(\tau_1, \tau_2) \in K_v$, if $\tau_1$ and $\tau_2$ have a common vertex;  

$(\tau_1, \tau_2) \in K_r$, if $\tau_1$ and $\tau_2$ are equal-residue;  

$(\tau_1, \tau_2) \in K_m$, if $\tau_1$ and $\tau_2$ are equal-modulus.

Note that the introduced subsets may intersect.

\begin{example} 1) Consider the pair  $(\tau_1, \tau_2)$, where $\tau_1 = \tau_{0(2), 1(2)}$, $\tau_2 = \tau_{0(2), 1(4)}$. Then  $(\tau_1, \tau_2) \in K_v \cup K_r$; that is, the class transpositions have a common vertex and are equal‑residue.

2) Consider the pair $(\tau_1, \tau_3)$, where $\tau_1 = \tau_{0(2), 1(2)}$, $\tau_3 = \tau_{0(3), 1(3)}$. Then $(\tau_1, \tau_3) \in K_g \cup K_r$; that is, both class transpositions are horizontal and equal‑residue.

3) Consider the pair $(\tau_4, \tau_5)$, where $\tau_4 = \tau_{1(4), 2(6)}$, $\tau_5 = \tau_{1(6), 2(4)}$. Then  $(\tau_1, \tau_3) \in K_r \cup K_m$; that is, the class transpositions are simultaneously equal‑residue and equal‑modulus.
\end{example}

\begin{theorem}  \label{4.3}
1) Let $\tau_{l(m_1),r(m_2)}$ and $\tau_{l(m_3),r(m_4)}$ be equal‑residue class transpositions and suppose that $\frac{m_1}{m_2}\neq \frac{m_3}{m_4}$. Then $ord(\tau_{l(m_1),r(m_2)}\cdot \tau_{l(m_3),r(m_4)})=\infty$.

2) Let $\tau_{r_1(m),r_2(n)}$ and $\tau_{r_3(m),r_4(n)}$ be equal‑modulus class transpositions and suppose that $r_1\leqslant r_4,\,r_3\leqslant r_2$. Then $ord(\tau_{r_1(m),r_2(n)}\cdot\tau_{r_3(m),r_4(n)})\in \{1,2,3,6,\infty\}$.
\end{theorem}

\begin{proof}
1) We will show that the graph $\Gamma(\tau_{l(m_1),r(m_2)},\tau_{l(m_3),r(m_4)})$  contains connected components of length greater than $n$ for any $n\in\mathbb{N}$.  Consider the system of Diophantine equations  

\begin{equation*}
\begin{cases}
l+m_1x_1=l+m_3x_2,\\
r+m_4x_2=r+m_2x_3,\\
l+m_1x_3=l+m_3x_4,\\
\cdots\cdots\cdots\cdots\cdots\cdots\cdots\\
r+m_4x_{2n}=r+m_2x_{2n+1}.
\end{cases}
\end{equation*}

If this system is solvable and has a solution in which $x_{2s}\neq x_{2t}$ for $s\neq t$, then the graph $\Gamma(\tau_{l(m_1),r(m_2)},\tau_{l(m_3),r(m_4)})$ contains the subgraph
$$
b_{x_1}a_{x_1}c_{x_2}d_{x_2}b_{x_3}a_{x_3}c_{x_4}d_{x_4}\cdots b_{x_{2n-1}}a_{x_{2n-1}}c_{x_{2n}}d_{x_{2n}}b_{x_{2n+1}}
$$
for any natural number $n$.

It is not difficult to verify that the system has a solution 
$$
(x_1, x_2, \ldots, x_{2n+1}) = 
$$
$$
= \left(x_1,\frac{m_1}{m_3}\cdot x_1,\ldots,\left(\frac{m_4}{m_2}\right)^{k-1}\left(\frac{m_1}{m_3}\right)^k\cdot x_1,\left(\frac{m_4}{m_2}\right)^k\left(\frac{m_1}{m_3}\right)^k\cdot x_1,\ldots,\left(\frac{m_4}{m_2}\right)^n\left(\frac{m_1}{m_3}\right)^n\cdot x_1\right),
$$ 
where $x_1=(m_2m_3)^nh$, $h\in\mathbb{Z}$.  

We will show that $x_{2s}\neq x_{2t}$ for $s\neq t$. Suppose, to the contrary, that 
$$
\left(\frac{m_4}{m_2}\right)^{s-1}\left(\frac{m_1}{m_3}\right)^s=\left(\frac{m_4}{m_2}\right)^{k-1}\left(\frac{m_1}{m_3}\right)^k
$$
Then we obtain
$$
1=\left( \frac{m_1m_4}{m_2m_3}\right)^{k-s},
$$
which is possible if and only if $m_1m_4=m_2m_3$. This contradicts the hypothesis. Similarly, one can show that $x_{2s+1}\neq x_{2t+1}$ for $s\neq t$.

Thus, we have shown that the graph $\Gamma(\tau_{l(m_1),r(m_2)},\tau_{l(m_3),r(m_4)})$ contains a connected component of length at least $\geq 2n$  for any $n\in\mathbb{N}$. According to  Theorem \ref{t1}, the product $\tau_{l(m_1),r(m_2)}\cdot \tau_{l(m_3),r(m_4)}$ has infinite order.

2) If $r_1=r_3$ or $r_2=r_4$, then we obtain the product of two class transpositions with a common vertex. By Theorem \ref{3.1}, we have $ord(\tau_{r_1(m),r_2(n)}\cdot\tau_{r_3(m),r_4(n)})\in \{1,3,\infty\}$.

Now assume that $r_1\neq r_3$ and $r_2\neq r_4$. Then $r_1(m)\cap r_3(m)=r_2(n)\cap r_4(n)=\emptyset$. If, in addition, either $r_1(m)\cap r_4(n)=\emptyset$ or $r_2(n)\cap r_3(m)=\emptyset$, then $ord(\tau_{r_1(m),r_2(n)}\cdot\tau_{r_3(m),r_4(n)})\in \{2,3,6\}$. Suppose $r_1(m)\cap r_4(n)\neq\emptyset$ and $r_2(n)\cap r_3(m)\neq\emptyset$. If $m=n$,  then we obtain the product of two horizontal class transpositions. It is easy to see that in this case $ord(\tau_{r_1(m),r_2(n)}\cdot\tau_{r_3(m),r_4(n)})\in \{2,3\}$.

Assume now that $m<n$. We construct the graph 
$$
\Gamma(\tau_{r_1(m),r_2(n)},\tau_{r_3(m),r_4(n)})=(V,E)
$$ 
and describe the connected components that may arise. By definition, $\{a_{k_i},b_{k_i}\}\in E$ and $\{c_{k_j},d_{k_j}\}\in E$ for any $i,\,j\in\mathbb{Z}$. Since 
$$
r_1(m)\cap r_3(m)=r_2(n)\cap r_4(n)=\emptyset,~~r_1(m)\cap r_4(n)\neq\emptyset~~\mbox{and}~~r_2(n)\cap r_3(m)\neq\emptyset,
$$ 
then between the graphs corresponding to the edges $a_{k_i}b_{k_i}$ and $c_{k_j}d_{k_j}$ either there are no edges at all, or $\{a_{k_i},d_{k_j}\}\in E$, or $\{b_{k_i},c_{k_j}\}\in E$. Therefore, the lengths of the cycles contained in the product $\tau_{r_1(m),r_2(n)}\cdot\tau_{r_3(m),r_4(n)}$, depend on the existence of solutions to systems of the following form:

\begin{equation*}
 \begin{cases}
   a_{x_1}=d_{x_2},\\
   c_{x_2}=b_{x_3},\\
   a_{x_3}=d_{x_4},\\
   \cdots\cdots\cdots\\
   c_{x_k}=b_{x_{k+1}},
 \end{cases} 
 \begin{cases}
   a_{x_1}=d_{x_2}, 
   \\
   c_{x_2}=b_{x_3},
   \\
   a_{x_3}=d_{x_4},
   \\
   \cdots\cdots\cdots
   \\
   a_{x_k}=d_{x_{k+1}},
 \end{cases}
\begin{cases}
   c_{x_1}=b_{x_2}, 
   \\
   a_{x_2}=d_{x_3},
   \\
   c_{x_3}=b_{x_4},
   \\
   \cdots\cdots\cdots
   \\
   c_{x_k}=b_{x_{k+1}},
 \end{cases} 
  \begin{cases}
   c_{x_1}=b_{x_2}, 
   \\
   a_{x_2}=d_{x_3},
   \\
   c_{x_3}=b_{x_4},
   \\
   \cdots\cdots\cdots
   \\
   a_{x_k}=d_{x_{k+1}}.
 \end{cases}
\end{equation*}
We consider the first of these systems (the other three are handled analogously). Rewriting the system, we have:
\begin{equation*}
\begin{cases}
mx_1-nx_2=r_4-r_1,\\
mx_2-nx_3=r_2-r_3,\\
mx_3-nx_4=r_4-r_1,\\
\cdots\cdots\cdots\cdots\cdots\cdots\cdots\\
mx_k-nx_{k+1}=r_2-r_3.
\end{cases}
\end{equation*}
Divide each equation by  $gcd(m,n)$ and set 
$$
\tilde{m}=\frac{m}{gcd(m,n)},~~\tilde{n}=\frac{n}{gcd(m,n)}.
$$
Define
$$
c_i=\frac{r_4-r_1}{gcd(m,n)}~~\mbox{for even}~~i,~~c_j=\frac{r_2-r_4}{gcd(m,n)}~~\mbox{for odd}~~j.
$$ 
Then the system becomes:
\begin{equation*}
\begin{cases}
\tilde{m}x_1-\tilde{n}x_2=c_1,\\
\tilde{m}x_2-\tilde{n}x_3=c_2,\\
\tilde{m}x_3-\tilde{n}x_4=c_3,\\
\cdots\cdots\cdots\cdots\cdots\\
\tilde{m}x_k-\tilde{n}x_{k+1}=c_k.
\end{cases}
\end{equation*}

We now show that this system has a solution for any $k$. Consider the first equation. It has a solution of the form 
$$
x_1(t_1)=\tilde{x}_1+\tilde{n}t_1, ~~x_2(t_1)=\tilde{x}_2+\tilde{m}t_1,
$$ 
where $\tilde{x}_1, \tilde{x}_2$ are particular solutions and $t_1\in\mathbb{Z}$. Substituting $x_2(t_1)$ into the second equation yields 
$$
\tilde{m}^2t_1-\tilde{n}x_3=c_2-\tilde{m}\tilde{x}_2.
$$
Since $gcd(\tilde{m}^2,\tilde{n})=gcd(\tilde{m},\tilde{n})=1$,   this equation has solutions 
$$
t_1(t_2)=\tilde{t}_1+\tilde{n},x_3(t_2)=\tilde{x}_3+\tilde{m}^2t_2,
$$ 
where $\tilde{t}_1, \tilde{x}_3$ are particular solutions and $t_2\in\mathbb{Z}$. Then the system
\begin{equation*}
 \begin{cases}
   \tilde{m}x_1-\tilde{n}x_2=c_1, 
   \\
   \tilde{m}x_2-\tilde{n}x_3=c_2,
 \end{cases}
\end{equation*}
has solutions
$$
(x_1, x_2, x_3) = (x_1(\tilde{t}_1+\tilde{n}t_2), x_2(\tilde{t}_1+\tilde{n}t_2), x_3(t_2)), ~~t_2\in\mathbb{Z}.
$$
Continuing this process, at the $l-$th step substituting $x_{l+1}(t_{l})$ into the equation $\tilde{m}x_{l+1}-\tilde{n}x_{l+2}=c_s$, we obtain a solution for the system consisting of $n$ equations.

Next, we show that $x_i\neq x_j$ if $i$ and $j$ have the same parity. This implies that all subgraphs of the form $a_{x_i}b_{x_i}$ and $a_{x_j}b_{x_j}$(and similarly, all subgraphs of the form $c_{x_i}d_{x_i}$ and $c_{x_j}d_{x_j}$) are distinct when $x_i, x_j$ are distinct solutions of the system.
Suppose to the contrary that there exist indices $i$ and $j$ of the same parity such that $x_i=x_j$. Then, considering the system 
\begin{equation*}
\begin{cases}
\tilde{m}x_1-\tilde{n}x_2=c_1,\\
\cdots\cdots\cdots\cdots\cdots\\
\tilde{m}x_i-\tilde{n}x_{i+1}=c_i,\\
\cdots\cdots\cdots\cdots\cdots\\
\tilde{m}x_{j-1}-\tilde{n}x_j=c_{j-1},\\
\cdots\cdots\cdots\cdots\cdots\\
\tilde{m}x_k-\tilde{n}x_{k+1}=c_k,
\end{cases}
\end{equation*}
we have in particular
\begin{equation*}
\begin{cases}
\tilde{m}x_i-\tilde{n}x_{i+1}=c_i,\\
\cdots\cdots\cdots\cdots\cdots\\
\tilde{m}x_{j-1}-\tilde{n}x_j=c_{j-1}.\\
\end{cases}
\end{equation*}
Since $x_i=x_j$, this becomes
\begin{equation*}
\begin{cases}
\tilde{m}x_i-\tilde{n}x_{i+1}=c_i,\\
\cdots\cdots\cdots\cdots\cdots\\
\tilde{m}x_{j-1}-\tilde{n}x_i=c_{j-1}.\\
\end{cases}
\end{equation*}
Note that the $x_t$ are either all positive or all negative (this follows from the definition of class transpositions). Without loss of generality, assume that all $x_t$ are positive. Then  
\begin{equation*}
\begin{cases}
\tilde{m}x_i>\tilde{n}x_{i+1},\\
\cdots\cdots\cdots\cdots\\
\tilde{m}x_{j-1}>\tilde{n}x_i.
\end{cases}
\end{equation*}
Since $\tilde{m}<\tilde{n}$, it follows that
\begin{equation*}
\begin{cases}
x_i>x_{i+1},\\
\cdots\cdots\cdots\cdots\\
x_{j-1}>x_i.
\end{cases}
\end{equation*}
This is a contradiction.

According to Theorem \ref{t1}, the product of the class transpositions contains cycles of length at least  $\geq k$ for every natural number $k$, and hence has infinite order.
\end{proof}

\bigskip

%%%%%%%%%%%%%%%%%%%%%%%%%%%%%%%%%%
\section{Acknowledgement}
This work was supported by a grant from the Russian Science Foundation, RSF
 24-11-00119. We are grateful to S. Kohl for his comments on a draft version of our manuscript and for his helpful suggestions.

\bigskip

%%%%%%%%%%%%%%%%%%%%%%%%%%%%%%%%%%

\end{document}